\documentclass[12pt]{article}%
\usepackage{amsmath,amssymb,amsthm,amsfonts}
\usepackage{wasysym}
\usepackage{graphicx}
\usepackage[dvipsnames]{xcolor}
\usepackage{stackengine}

\usepackage[colorlinks]{hyperref}
\usepackage{tikz}
\usepackage[export]{adjustbox}

\usepackage{appendix}

\newcommand{\ea}{\textit{et al. }}
\renewcommand{\epsilon}{\varepsilon}

\definecolor{red}{rgb}{0.8500, 0.3250, 0.0980}
\definecolor{green}{rgb}{0.4660, 0.6740, 0.1880}
\definecolor{yellow}{rgb}{0.9290, 0.6940, 0.1250}
\definecolor{blue}{rgb}{0, 0.4470, 0.7410}

\begin{document}

\title{Stochastic discrete dynamical model for the hydrodynamic analog of a quantum mirage}

\author{Gonzalo Ferrandez Quinto \thanks{Department of Applied Mathematics, University of Washington} \thanks{Department of Physics, University of Washington}, Aminur Rahman\footnotemark[1] \thanks{Corresponding  \url{arahman2@uw.edu}}}

\date{}

\maketitle

\begin{abstract}
Recently discrete dynamical models of walking droplet dynamics have allowed for fast numerical simulations of their horizontal chaotic motion and consequently the long-time statistical distribution of the droplet position.  We develop a discrete model for walkers on an elliptical corral with two dominant eigenmodes.  These eigenmodes are excited periodically with each having random, but correlated, weights.  At each iteration the model computes the horizontal propulsion of the droplet due to impacts with the fluid bath.  The propulsion is calculated based on the interaction of the droplet with the wavefield produced by the eigenmode excitations.  We record the long-time statistics, and incredibly for such a simple model, they are qualitatively similar to that of experiments.
\end{abstract}

\section{Introduction}
\label{Sec: Intro}



Walking droplets have captured the imagination of classical physicists over the past two decades.  It has been known since the 70s that, with the right conditions, a liquid droplet can bounce for long periods on a vertically vibrating liquid bath \cite{Walker1978}.  Later, Couder and coworkers showed the droplet can start moving horizontally, and thus dubbed ``walkers'', if the forcing on the fluid bath is increased \cite{CPFB05, CouderFort06}.  Since a particle (the droplet) is being propelled by the waves it creates on the bath, researchers have been able to produce quantum-like phenomena such as a particle confined in a circular corral \cite{CLE1993, HMFCB13, CSB18} and the quantum mirage \cite{QuantumMirage, Saenz-Mirage2018}, among others that have been well delineated in reviews by Bush and coworkers \cite{Bush10, Bush15a, Bush15b, BushOza20_ROPP}.
\begin{figure}[htbp]
\centering
\includegraphics[width=0.9\textwidth]{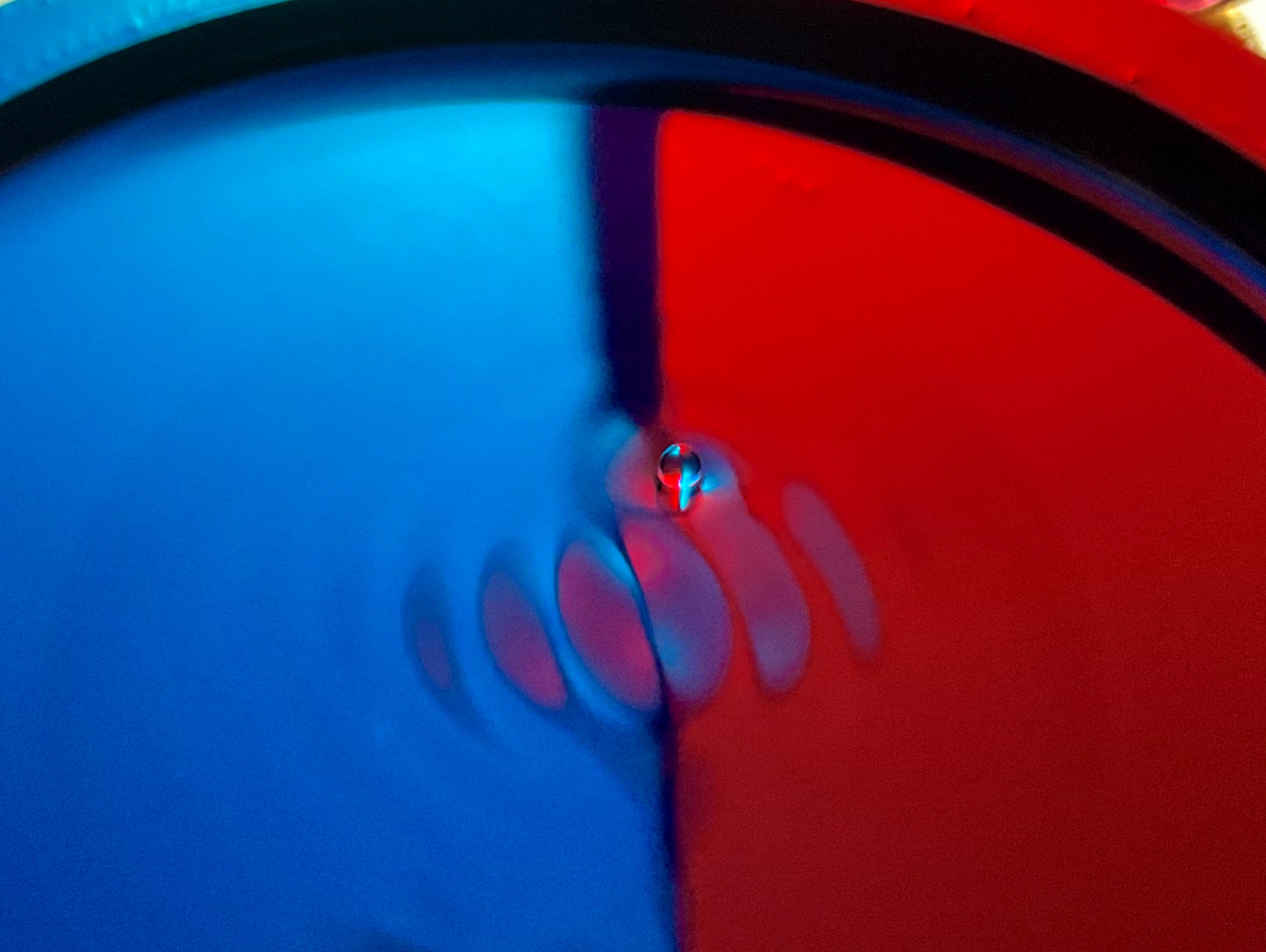}
\caption{A walking droplet producing waves on a vibrating fluid bath. Blue and red light filters are present to improve the visualization of the wave field.}
\label{Fig: droplet}
\end{figure}
%


While the phenomena admitted by the droplets system is intriguing, the modeling has proven to be quite delicate and overwhelmingly complex.  Even reduced continuum models, such as stroboscopic models \cite{ORB13}, can be computationally expensive and mathematically intractable.  Discrete dynamical models \cite{Gilet14, RahmanBlackmore16, Gilet16, Rahman18} provide a respite from complexity and allow for rigorous mathematical analysis \cite{RJB17, RahmanBlackmore20, RahmanKutzDampedDriven}.  The first such model, by Gilet \cite{Gilet14}, confines a walker in a 1-dimensional domain.  The wave and position of the droplet are updated iteratively at each impact with interaction terms between the two.  This enabled Gilet to conjecture \cite{Gilet14} and Rahman and Blackmore to prove \cite{RahmanBlackmore16} the existence of Neimark--Sacker bifurcations \cite{Neimark, Sacker}.  Later models of this type revealed new heteroclinic-type bifurcations and routes to chaos \cite{RJB17, RahmanBlackmore20}.  Gilet also developed a discrete dynamical model for a walker confined in a ciruclar corral \cite{Gilet16}, which incredibly reproduced the longtime statistics shown in the hydrodynamic corral experiments \cite{HMFCB13}.  This communication aims to modify Gilet's circular corral model for the hydrodynamic analog of the quantum mirage that was recently developed by S\'{a}enz \textit{et al.} \cite{Saenz-Mirage2018}.

In this experiment, an elliptic corral is vibrated vertically with a given amplitude of acceleration $\gamma$. When the driving acceleration is below the Faraday wave threshold, $\gamma/\gamma_F<1$, standing waves will only be visible under the presence of a droplet. In order to capture which modes could be excited by the walker, S\'{a}enz \textit{et al.} \cite{Saenz-Mirage2018} drove an elliptical corral of major semi-axis of 14.25mm and eccentricity of 0.5 at accelerations just above the Faraday wave threshold. This demonstrated that two distinct Mathieu function standing waves are preferentially excited in the sub Faraday regime. These functions correspond to the Helmholtz equation solutions in an elliptical domain with fixed boundary conditions. Specifically the even mode (n=4, j=4) and odd mode (n=1, j=5) as described by Wilson \ea \cite{Wilson-MathieuEqs2007}.

Below the Faraday wave threshold, the walking droplet experiment was performed with a uniform corral depth. For this bottom topography, an equally weighted superposition of both dominant modes resembled the average wave pattern recorded at the surface of the corral. To understand the effects of depth variability on the wave field, the experiment was repeated with an increased depth impurity at one ellipse focus. It was observed that this resulted in a wave field heavily dominated by the even mode (4,4). Further, when the depth impurity was located in the middle of the semi-minor axis, it resulted in a reduction in the resonance of both modes and therefore a lower general velocity of the walker. In addition, the horizontal characteristic pattern of the odd mode (1,5) was reinforced in the observed wave field. In both cases, the impurity was observed to act as attraction point for the walker.

The remainder of the paper is as follows:  in Sec. \ref{Sec: Dynamics} we model the experiment of S\'{a}enz \textit{et al.} as a discrete dynamical system.  Then, Sec. \ref{Sec: Numerics} delineates the numerical methods used to simulate and analyze the model.  In Sec. \ref{Sec: Results} we present results comparing the model to experiments.  Finally, in Sec. \ref{Sec: Conclusion} we conclude with some closing remarks and future directions.


\section{Discrete dynamical model}
\label{Sec: Dynamics}

In order to model the hydrodynamic analog of the quantum mirage \cite{Saenz-Mirage2018} as a discrete dynamical system, we must first study the models of Gilet \cite{Gilet14, RahmanBlackmore16, Gilet16}.  The 1-dimensional model \cite{Gilet14}
\begin{subequations}
    \begin{align}
        w_{n+1} &= \mu[w_n + \Psi(x_n)],\\
        x_{n+1} &= x_n - Cw_n\Psi'(x_n);
    \end{align}
    \label{Eq: Gilet 1-D}
\end{subequations}
describes the discrete evolution at each impact $n$ of the wavefield with damping factor $\mu \in [0,1]$, where $w$ is the amplitude and $\Psi$ is the shape of the eigenmode.  Further, $x$ is the horizontal position and $C \in [0,1]$ is the wave-particle coupling factor, which is related to the coefficient of restitution of the inelastic collisions.  In the model it is assumed that the droplet moves opposite to the gradient of the wavefield at impact, which is heuristically reasonable if we consider a point mass impacting a rigid surface.  Although the droplet-surface interaction is more nuanced for the fluidic system, in 1-D, opposite the gradient is certainly a reasonable assumption due to the low kinematic degree of freedom.

This idea was extended to the circular corral by Gilet \cite{Gilet16} with the model
\begin{subequations}
    \begin{align}
        r_{n+1} &= \frac{r_n - C\sum_k \frac{d\varphi}{dr}\bigg|_{r_n}\text{Re}\left(w_{k,n}e^{ik\theta_n}\right)}{\cos(\theta_{n+1} - \theta_n)},\\
        r_{n+1} &= \frac{\frac{C}{r_n}\sum_k \varphi'(r_n)\text{Im}\left(w_{k,n}e^{ik\theta_n}\right)}{\sin(\theta_{n+1} - \theta_n)},\\
        w_{k,n+1} &= \mu_k\left[w_{k,n} + \phi_k(r_n)e^{-ik\theta_n}\right];
    \end{align}
    \label{Eq: Gilet 2-D}
\end{subequations}
where, in addition to the setup of \eqref{Eq: Gilet 1-D}, we now have to keep track of multiple eigenmodes that can be excited in 2-dimensions, which are represented by the integer $k$, and the value of $C$ is no longer restricted to be less than unity.  We also use the standard $(r,\theta)$ polar coordinate system.  Finally, $\varphi$ represent the eigenmodes that can be excited on a disk.  In this communication we modify \eqref{Eq: Gilet 2-D} to model the walker on an ellipse with an impurity.

%
We first solve for the possible eigenmodes.  We begin by writing Helmholtz equation,
\begin{equation}
\nabla^{2}\Psi=-k^2\Psi,
\end{equation}
which arises when solving the wave equation by separation of variables.  In this case $k$ is a separation constant and physically represents the angular wavenumber. Solving the Helmholtz equation on a disk yields Bessel functions. Solving it on an ellipse with the elliptical coordinates yields Mathieu functions. To convert from elliptical coordinates to Cartesian we use
\begin{align*} 
x&=A\ \cosh \xi \ \cos \eta \\
y&=A\ \sinh \xi \ \sin \eta
\end{align*}
where $A = 7.125\text{mm}$ from experiments \cite{Saenz-Mirage2018}, and boundary
\begin{equation*}
0\leq\xi\leq\xi_{0},\qquad -\pi\leq\eta\leq\pi,
\end{equation*}
where $\xi_{0}=\tanh^{-1}(b/a)$ with $a$ and $b$ representing the major and minor semi-axes values of the elliptical corral. When the Helmholtz equation is solved using separation of variables in elliptical coordinates, the equation depending on $\xi$ is called the radial Mathieu equation, and the equation depending on $\eta$ is the angular Mathieu equation. As shown in (6) through (9) in the work of Wilson \ea \cite{Wilson-MathieuEqs2007}, separation of variables introduces constants $n$ and $q$. Due to periodic boundary conditions of $\eta$, the value of $n$ must be a whole number. Solutions for the angular equation can be found through a recurrence relation. For the radial Mathieu function, solutions can be found through an infinite expansion of Bessel functions of the first kind. The no-slip boundary condition on the corral introduce further constraints on the variable $q$. A value of $q$, for a given value of $n$, must satisfy $\Psi_{n,j}(\xi_0,\eta,q_{n,j})=0$, where $q_{n,j}$ is the $j^{th}$ value for which this condition is satisfied. Further, functions of the form $\Psi_{n,j}(\xi,\eta,q_{n,j})=\Phi(\xi,q_{n,j})\Theta(\eta,q_{n,j})$ have even symmetry when $\Theta(\eta,q_{n,j})$ is expanded with a finite series of cosines and odd symmetry when expanded with sine functions.

As observed by S\'{a}enz \ea \cite{Saenz-Mirage2018}, the modes are
\begin{subequations}
    \begin{align}
    \Psi_{1,5}(\xi,\eta,q_{1,5})=Ms_1^{(1)}(\xi,q_{1,5})se_1^{(1)}(\eta,q_{1,5})\\
    \Psi_{4,4}(\xi,\eta,q_{4,4})=Mc_4^{(1)}(\xi,q_{4,4})ce_4^{(1)}(\eta,q_{4,4})
    \end{align}
    \label{Eq: modes}
\end{subequations}
\begin{figure}[htbp]
\centering
\stackinset{l}{}{t}{1pt}{\textbf{(a)}}{\includegraphics[width=.47\linewidth]{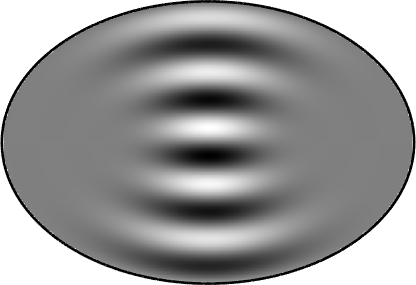}}
\stackinset{l}{}{t}{1pt}{\textbf{(b)}}{\includegraphics[width=.47\linewidth]{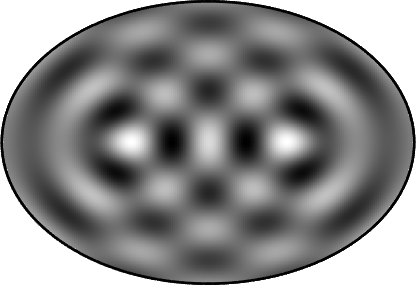}}
\caption{The Mathieu function resonant modes for an ellipse with semi-major axis of 14.25mm and eccentricity of 0.5. Lighter regions represent more positive mode values while darker regions represent more negative mode values. \textbf{(a)} Profile of the odd mode $\Psi_{1,5}(\xi,\eta,q_{1,5})$ where $n=1$, $j=5$ and $q_{1,5}$ is evaluated numerically to be $q_{1,5}=20.6638$. \textbf{(b)} Profile of the even mode $\Psi_{4,4}(\xi,\eta,q_{4,4})$ were $n=4$, $j=4$ and $q_{4,4}$ is evaluated numerically to be $q_{4,4}=21.4267$.}
\label{Fig: Pure Mathieu modes}
\end{figure}
The first odd mode is numerically evaluated to be $q_{1,5}=20.6638$, and the second even mode is $q_{4,4}=21.4267$. Examples of these eigenmodes excited by the Mathieu equation are shown in Fig \ref{Fig: Pure Mathieu modes}. In the heatmap, the dark regions represent negative values and light regions represent positive values. The domain is an ellipse with major semi-axis of $14.25\text{mm}$ and eccentricity of $0.5$.  The odd mode is plotted in Fig. \ref{Fig: Pure Mathieu modes}(a) and the even mode is plotted in Fig. \ref{Fig: Pure Mathieu modes}(b).


Using Gilet's model \cite{Gilet16} as inspiration let us first write the droplet motion and wavefield evolution with the same basic structure as \eqref{Eq: Gilet 2-D} in Cartesian coordinates,
    \begin{align*}
        x_{n+1} &= x_n - Cw_{n+1}\Psi_x(x_n,y_n),\\
        y_{n+1} &= y_n - Cw_{n+1}\Psi_y(x_n,y_n),\\
        w_{n+1} & = \mu\left(w_n + \Psi(x_n,y_n)\right).
    \end{align*}
We recall that this model assumes the droplet moves opposite the gradient of the wavefield.  While this is a reasonable assumption, when we simulate the model we observe that the long-time statistics are not close to that of experiments of S\'{a}enz \textit{et al.} \cite{Saenz-Mirage2018}.  In fact, there has yet to be a successful model for the mirage experiment, all of which make the opposite the gradient assumption.  There are several possibilities for these discrepancies:  perhaps the dominant eigenmodes observed in experiments are in actuality slightly different, the opposite wavefield gradient assumption may not hold in all cases, or a completely new modeling paradigm is necessary.

\begin{figure}[htbp]\centering
\stackinset{l}{}{t}{1pt}{\textbf{(a)}}{\includegraphics[width = 0.3\textheight]{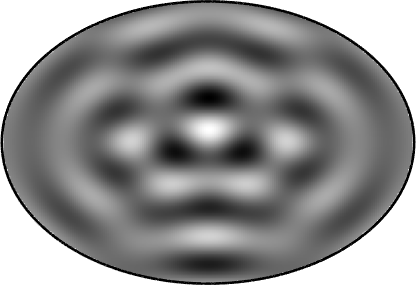}}
\stackinset{l}{}{t}{1pt}{\textbf{(b)}}{\includegraphics[width = 0.3\textheight]{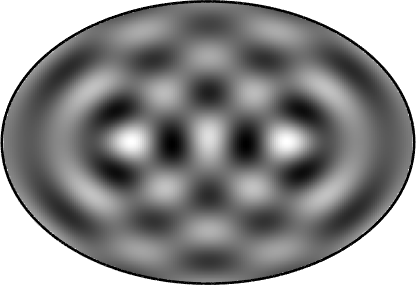}}
\stackinset{l}{}{t}{1pt}{\textbf{(c)}}{\includegraphics[width = 0.3\textheight]{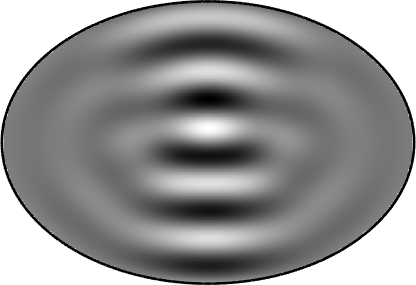}}
\caption{Normalized time averaged wavefields from simulations on the elliptical corral following \eqref{Eq: Wavefield}. Lighter regions represent more positive wavefield values while darker regions represent more negative wavefield values. Simulated wavefield for \textbf{(a)} has equally weighted eigenmodes $\alpha = 0.5$ and $\beta=0.5$ in \eqref{Eq: Wavefield},  \textbf{(b)} the $\Psi_{4,4}$ eigenmode dominates with $\alpha = 0.05$ and $\beta=0.5$ in \eqref{Eq: Wavefield}, and \textbf{(c)} the $\Psi_{1,5}$ eigenmode dominates with $\alpha = 0.5$ and $\beta=0.1$ in \eqref{Eq: Wavefield}.}
\label{Fig: Simulated average wavefields}
\end{figure}

In this work we hypothesize the droplet moves not in the direction opposite the gradient of the wavefield at impact, but rather in the perpendicular direction.  Mathematically this yields the model
\begin{subequations}
    \begin{align}
        x_{n+1} &= x_n - Cw_{n+1}\Psi_y(x_n,y_n),\\
        y_{n+1} &= y_n + Cw_{n+1}\Psi_x(x_n,y_n),\\
        w_{n+1} & = \mu\left(w_n + \Psi(x_n,y_n)\right);
    \end{align}
    \label{Eq: Our Model}
\end{subequations}
where $\Psi$ is a linear combination of two dominant eigenmodes from the solution to the Mathieu equation as observed by Saenz \textit{et al.} \cite{Saenz-Mirage2018}.  While the precise form of the linear combination was not observed in experiments, we assume it takes the form
\begin{equation}
    \Psi(x_n,y_n) = p_n\alpha\Psi_{1,5}+(1/2-p_n)\beta\Psi_{4,4}
    \label{Eq: Wavefield}
\end{equation}
In this case $p_n$ is a uniformly distributed random variable in the range $[0,0.5]$. This random variable is recalculated at each update of the droplet location.  The constants $\alpha$ and $\beta$ are the static weights assigned to each one of the modes. The relative variation of these parameters is responsible for producing the different time averaged wavefields of Fig. \ref{Fig: Simulated average wavefields}.
%
%


\section{Numerical methods}
\label{Sec: Numerics}

To produce the odd and even angular Mathieu functions we use the MATLAB toolbox \cite{matlab} associated with the article of Wilson \ea  \cite{Wilson-MathieuEqs2007}.  In this case, $ce_4^{(1)}$ and $se_1^{(1)}$, the angular Mathieu functions can be expanded using an infinite series of sine or cosine functions. Assuming that the series converges, a tri-diagonal matrix eigenvalue problem can be used to solve for the value of the coefficients as shown by Wilson \ea \cite{Wilson-MathieuEqs2007} to produce the odd and even radial Mathieu functions, $Ms^{(1)}_1(\xi,q_{1,5})$ and $ Mc_4^{(1)}(\xi,q_{4,4}) $. The terms are expanded using an infinite series of Bessel functions of the first kind \cite{Wilson-MathieuEqs2007}. The approximation uses the first fifty terms of the infinite series. There is no need to recalculate the coefficients of the Bessel series due to the relation of sinusoidal functions with hyperbolic functions.  Moreover, the toolbox employed for the Mathieu functions \cite{matlab} uses elliptical coordinates. In order to simplify our simulation and model, we only employ Cartesian coordinates and convert to elliptical coordinates when calculating the eigenmode values through the toolbox.


%
Now let us embed a grid with $h = \Delta x$ and $k = \Delta y$ onto our eigenmodes.  Since we wish to approximate the gradient at each timestep, and do not reuse it in subsequent approximations, we are not as concerned with the propagation of truncation errors.  Therefore, we approximate the gradient using second order centered differences,
\begin{align*}
    \frac{\partial \Psi}{\partial x}\bigg|_{(x_n, y_n)} &\approx \frac{\Psi(x_n + h, y_n) - \Psi(x_n - h, y_n)}{2h},\\
    \frac{\partial \Psi}{\partial y}\bigg|_{(x_n,y_n)} &\approx \frac{\Psi(x_n, y_n + k) - \Psi(x_n, y_n - k)}{2k}.
\end{align*}

To iterate the dynamical system, we first update the wavefield as a linear combination of the two dominant eigenmodes with a random, but correlated, amplitude in the form of \eqref{Eq: Wavefield}.  Then we calculate the gradient using the centered difference method.  Next we update the wavefield followed by the position according to \eqref{Eq: Our Model}.  Occasionally the droplet escapes the boundary since boundary effects are not taken into consideration. In these instances, the simulation is considered to be ended and a new run is started at another random point of the corral. The final histogram is the cumulative summation of several (on the order of ten) different runs without boundary effects. We iterate the model on the order of $10^5$ times in total, which gives us enough points to produce a heat map of the long-time statistics presented in the sequel.  The heat maps are produced by using $90 \times 90$ bins in a histogram of the droplet location in the interior of the ellipse.


\section{Results}
\label{Sec: Results}

In this section we simulate the wavefield and droplet motion, which is then compared to the experiments of S\'{a}enz \ea \cite{Saenz-Mirage2018}.  The wavefield is given by \eqref{Eq: Wavefield} and the droplet motion is modeled as \eqref{Eq: Our Model}.  The wavefields themselves are derived directly from the hypotheses of S\'{a}enz \ea \cite{Saenz-Mirage2018}, and will naturally agree with experiments.  The droplet motion, on the other hand, is quite a bit trickier, and we seek qualitative agreements with experiments.  Although we propel the droplet perpendicular to the gradient in our model \eqref{Eq: Our Model} rather than the standard opposite to the gradient, we still observe good qualitative agreement.

We first simulate the time averaged wavefields, given by \eqref{Eq: Wavefield}, in Fig. \ref{Fig: Simulated average wavefields}. Figure \ref{Fig: Simulated average wavefields}\textbf{(a)} is the case when there is no impurity in the bath, $\alpha=0.5$ and $\beta=0.5$. When the depth impurity is located at one of the foci, the values are $\alpha=0.05$ and $\beta=0.5$, corresponding to Fig. \ref{Fig: Simulated average wavefields}\textbf{(b)}. This follows the reasoning that, experimentally, the impurity reduces the excitation of the $\Psi_{1,5}$ mode due to the odd symmetry of the mode with respect to the x-axis \cite{Saenz-Mirage2018}. The $\Psi_{4,4}$ mode is unaffected as it has even symmetry with respect to the x-axis. When the impurity is at the center of the semi-minor axis, the pattern observed heavily resembles the $\Psi_{1,5}$ mode, which justifies the choice of $\alpha=0.5$ and $\beta=0.1$, corresponding to Fig. \ref{Fig: Simulated average wavefields}\textbf{(c)}.

\begin{figure}[htbp]\centering
\stackinset{l}{}{t}{1pt}{\textbf{(a)}}{\includegraphics[height = 0.22\textheight]{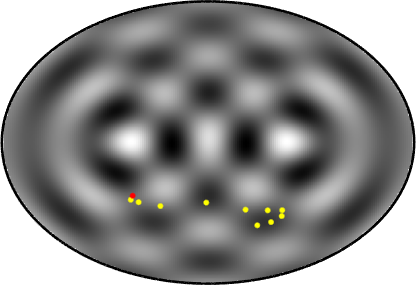}}
\stackinset{l}{}{t}{1pt}{\textbf{(b)}}{\includegraphics[height = 0.22\textheight]{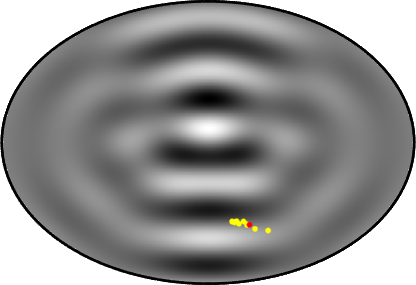}}
\stackinset{l}{}{t}{1pt}{\textbf{(c)}}{\includegraphics[height = 0.22\textheight]{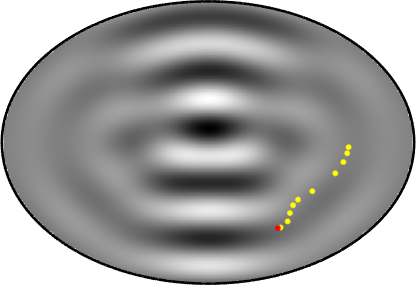}}
\stackinset{l}{}{t}{1pt}{\textbf{(d)}}{\includegraphics[height = 0.22\textheight]{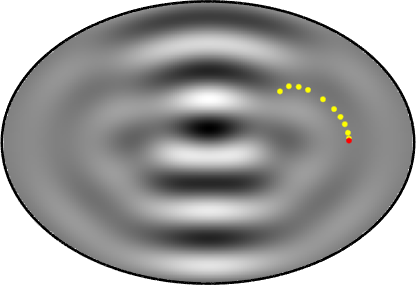}}

\caption{Simulated trajectory of the droplet on equally weighted eigenmodes, in the wavefield \eqref{Eq: Wavefield}, in the elliptical corral. Lighter regions represent more positive wavefield values while darker regions represent more negative wavefield values. The background plot corresponds to the effective instantaneous wavefield $w_{n+1}\Psi(x_n,y_n)$ at the initial iteration. The initial point corresponds to the red marker and the yellow markers correspond to a subset of  points of the followed path. figures \textbf{(a)} through \textbf{(d)} correspond to different time snapshots of a sample trajectory.}
\label{Fig: Simulated wavefields}
\end{figure}

The instantaneous wavefield at any given iteration can be interpreted as $w_{n+1}\Psi(x_n,y_n)$, with the horizontal trajectory governed by the model \eqref{Eq: Our Model} with wavefield \eqref{Eq: Wavefield} containing equally weighted parameters $\alpha = 0.5$ and $\beta=0.5$. A sample trajectory can be observed in Fig. \ref{Fig: Simulated wavefields} where the instantaneous wavefield of the bath is superimposed with several of the following iterations of the trajectory. The instantaneous wavefield corresponds to the red marker position of the droplet. A highly variable wavefield can be observed, which follows from the inclusion of the distributed random variable $p_{n}$ in \eqref{Eq: Wavefield}.  Remarkably, in Fig. \ref{Fig: Simulated wavefields}, the droplet follows a path that carves the local maximum of the instantaneous field in Figs. \ref{Fig: Simulated wavefields}\textbf{(c)} and \textbf{(d)}. This shows striking resemblance to trajectories present in Fig. 2\textbf{(a)} and \textbf{(b)} of S\'{a}enz \ea \cite{Saenz-Mirage2018}.
%

\begin{figure}[htbp]\centering
\stackinset{l}{}{t}{1pt}{\textbf{(a)}}{\includegraphics[height = 0.24\textheight]{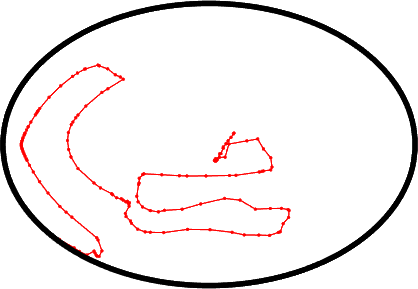}}
\stackinset{l}{}{t}{1pt}{\textbf{(b)}}{\includegraphics[height = 0.24\textheight]{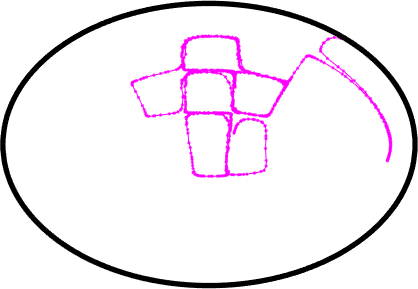}} 
\stackinset{l}{}{t}{1pt}{\textbf{(c)}}{\includegraphics[height = 0.24\textheight]{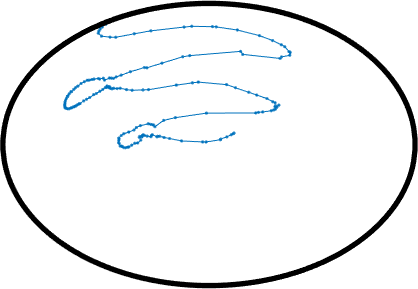}}
\caption{Examples of simulated droplet trajectories governed by \eqref{Eq: Our Model} on the elliptical corral.   \textbf{(a)}  Simulated droplet trajectories for equally weighted eigenmodes. $\alpha = 0.5$ and $\beta=0.5$ in \eqref{Eq: Wavefield}. \textbf{(b)} Simulated droplet trajectories for the $(4,4)$ dominant eigenmode with $\alpha = 0.05$ and $\beta=0.5$ in \eqref{Eq: Wavefield}.  \textbf{(c)} Simulated droplet trajectories for the $(1,5)$ dominant eigenmode with $\alpha = 0.5$ and $\beta=0.1$ in \eqref{Eq: Wavefield}.}
\label{Fig: Simulated trajectories}
\end{figure}
\begin{figure}[htbp]\centering
\stackinset{l}{}{t}{1pt}{\color{white}\textbf{(a)}}{\includegraphics[height = 0.24\textheight]{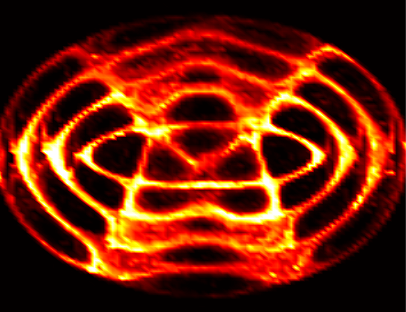}}\quad
\stackinset{l}{}{t}{1pt}{\color{white}\textbf{(b)}}{\includegraphics[height = 0.24\textheight]{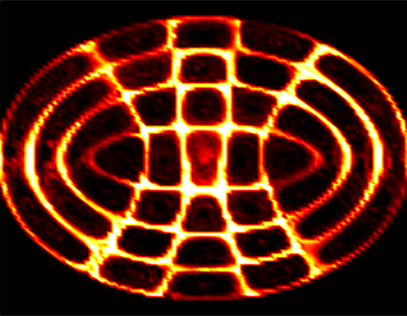}}\quad
\stackinset{l}{}{t}{1pt}{\color{white}\textbf{(c)}}{\includegraphics[height = 0.24\textheight]{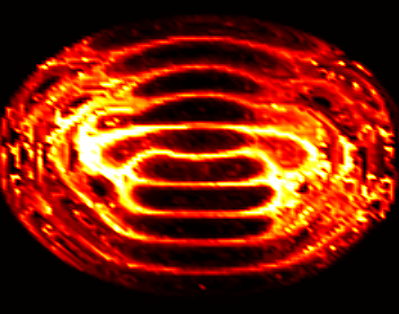}}
\stackinset{l}{}{t}{1pt}{}{\includegraphics[height = 0.24\textheight]{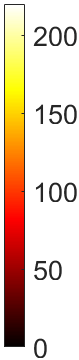}}
\caption{Long-time statistics of the simulated trajectories governed by \eqref{Eq: Our Model} on the elliptical corral.  Heat maps are produced by taking a histogram with $90 \times 90$ bins.  The top end of the colorbar (white) represents greater than $220$ occurrences for the relevant bins.  \textbf{(a)}  Long-time statistics for simulated droplet trajectories induced by equally weighted eigenmodes. $\alpha = 0.5$ and $\beta=0.5$ in \eqref{Eq: Wavefield}.  \textbf{(b)} Long-time statistics for simulated droplet trajectories induced by  the $(4,4)$ dominant eigenmode with  with $\alpha = 0.05$ and $\beta=0.5$ in \eqref{Eq: Wavefield}.  \textbf{(c)} Long-time statistics for simulated droplet trajectories induced by the $(1,5)$ dominant eigenmode with $\alpha = 0.5$ and $\beta=0.1$ in \eqref{Eq: Wavefield}.}
\label{Fig: Long-time statistics}
\end{figure}

Next we simulate the droplet motion, given by \eqref{Eq: Our Model}, and their associated long-time statistics in Figs. \ref{Fig: Simulated trajectories} and \ref{Fig: Long-time statistics}. Figure \ref{Fig: Simulated trajectories} corresponds to sample trajectories for different values of the weights in wavefield \eqref{Eq: Wavefield}. We observe that the qualitative behavior of the trajectories depend sensitively on the weights of the eigenmodes. In Figs. \ref{Fig: Simulated trajectories} \textbf{(a)} and \ref{Fig: Long-time statistics} \textbf{(a)}, the value of the weights correspond to the case where there is no depth impurity at the bottom of the bath. In that case, both eigenmodes $\Psi_{1,5}$ and $\Psi_{4,4}$ are considered to be excited equally with weights $\alpha=0.5$ and $\beta=0.5$. In Figs. \ref{Fig: Simulated trajectories} \textbf{(b)} and \ref{Fig: Long-time statistics} \textbf{(b)}, the value of the weights correspond to the case were there is a depth impurity at the bottom of the bath in one of the foci of the ellipse. In that case, eigenmode $\Psi_{1,5}$ is reduced and eigenmode $\Psi_{4,4}$ dominates with corresponding weights $\alpha=0.05$ and $\beta=0.5$. Finally, in Figs. \ref{Fig: Simulated trajectories} \textbf{(c)} and \ref{Fig: Long-time statistics} \textbf{(c)}, the value of the weights correspond to the case where there is a depth impurity at the middle of the semi-minor axis of the ellipse. In that case, eigenmode $\Psi_{1,5}$ dominates and eigenmode $\Psi_{4,4}$ is reduced with corresponding weights $\alpha=0.5$ and $\beta=0.1$. In all simulations both modes are needed in order for the droplet to be able to circulate around the elliptical bath. It is also worth noting that in Fig. \ref{Fig: Long-time statistics}, the long-time statistics are primarily dependant on the relative value of the parameters $\alpha$ and $\beta$.  We observe that the heat maps of the long-time statistics of Fig. \ref{Fig: Long-time statistics} are reminiscent of Figs. 2\textbf{(e)}, 5\textbf{(b)}, and 5\textbf{(a)} in S\'{a}enz \ea \cite{Saenz-Mirage2018}, respectively.
%

\begin{figure}[htbp]\centering
\stackinset{l}{}{t}{1pt}{\color{white}\textbf{(a)}}{\includegraphics[height = 0.27\textheight]{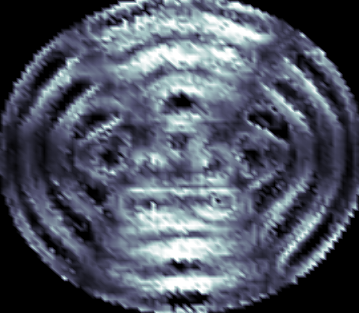}}\quad
\stackinset{l}{}{t}{1pt}{\color{white}\textbf{(b)}}{\includegraphics[height = 0.27\textheight]{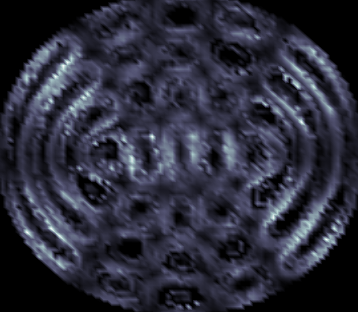}}\quad
\stackinset{l}{}{t}{1pt}{\color{white}\textbf{(c)}}{\includegraphics[height = 0.27\textheight]{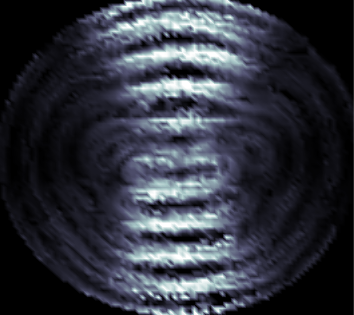}}
\stackinset{l}{}{t}{1pt}{}{\includegraphics[height = 0.27\textheight]{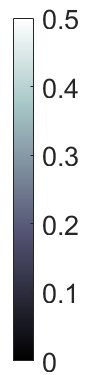}}
\caption{Long-time statistics of the average displacement in millimeters per iteration, governed by \eqref{Eq: Our Model} on the elliptical corral.  Heat maps are produced by taking a histogram with $90 \times 90$ bins. \textbf{(a)}  Long-time statistics of the average displacement in millimeters per iteration induced by equally weighted eigenmodes with $\alpha = 0.5$ and $\beta=0.5$ in \eqref{Eq: Wavefield}.  \textbf{(b)} Long-time statistics of the average displacement in millimeters per iteration induced by  the $(4,4)$ dominant eigenmode with $\alpha = 0.05$ and $\beta=0.5$ in \eqref{Eq: Wavefield}.  \textbf{(c)} Long-time statistics of the average displacement in millimeters per iteration induced by the $(1,5)$ dominant eigenmode with $\alpha = 0.5$ and $\beta=0.1$ in \eqref{Eq: Wavefield}.}
\label{Fig: Long-time-vel statistics}
\end{figure}
Similar to Fig. \ref{Fig: Long-time statistics}, Fig. \ref{Fig: Long-time-vel statistics} corresponds to the average displacement per iteration at different locations of the elliptical corral. These displacement histograms are derived from the position histograms in Fig. \ref{Fig: Long-time statistics}. Due to the lack of temporal dimension in the model, these graphs do not correspond exactly with velocity, but provides insight into the dynamics of the model. It is worth noting the inverse relation of Figs. \ref{Fig: Long-time statistics} and \ref{Fig: Long-time-vel statistics}, which is similar to what is seen in experiments \cite{Saenz-Mirage2018}. That is, locations where the position histogram shows a higher count tend to have a lower displacement per iteration value. 
%



\section{Discussion}
\label{Sec: Conclusion}


In this communication we model and simulate the hydrodynamic analog \cite{Saenz-Mirage2018} of the quantum mirage \cite{QuantumMirage}.  We first develop a stochastic discrete dynamical model \eqref{Eq: Our Model} inspired by the confined geometry \eqref{Eq: Gilet 1-D} \cite{Gilet14} and circular corral \eqref{Eq: Gilet 2-D} \cite{Gilet16} models of Gilet.  While Gilet's model assumes the droplet is propelled opposite the gradient between iterations, our model assumes that the droplet is propelled perpendicular to the gradient.  Indeed, instantaneously, we would expect the wavefield to exert a force in the opposite direction (in 2-D) to the gradient of the wavefield, however the iterated map need not represent each impact.  There could very well be multiple impacts between iterations.  Further, in the circular corral model of Gilet \eqref{Eq: Gilet 2-D} \cite{Gilet16}, several eigenmodes are used to model the wavefield. In our model two dominant eignemodes are excited with stochastically varying instantaneous weights that contribute to the wavefield at each timestep as given by \eqref{Eq: Wavefield}.  Although the numerics are considerably simplified by the use of a discrete dynamical system, there are some numerical considerations due to the use of the Mathieu equation.  Finally, we present simulations of the wavefield, droplet trajectories, and long-time statistics that qualitatively agree with the experiments of S\'{a}enz \ea \cite{Saenz-Mirage2018}.


In addition to the close qualitative agreement with experiments, the present model's appeal lies in the use of relative weights on the mode superposition to produce vastly different long-time statistics. Moreover, the assumption that the droplet movement is perpendicular to the gradient, ensures that the path followed by the simulated droplet shadows the time averaged wavefield as observed in the experiments of S\'{a}enz \ea \cite{Saenz-Mirage2018}.  Interestingly, even with this assumption, Fig. \ref{Fig: Long-time-vel statistics} shows that the displacement over one iteration, at least on average, is opposite the gradient.  Further, the model is wavefield agnostic; that is, we can excite wavefields in the model based on observations, and then let the model deterministically (over one iteration) propel the droplet.  Finally, the model is observed to be time reversal invariant.  Unlike more complex hydrodynamic models, we are not reliant on the infinite summation of wavefields to simulate memory in the system, and therefore can just as easily run our simulations backwards in time.  Hydrodynamically this may seem odd, but for certain quantum analogs this could help us gain insight into the connections with quantum systems.


Despite its many advantages, it would benefit from a rigorous fluid mechanical analysis of droplet motion in relation to the wavefield.  That is, why has opposite the gradient not worked in a variety of modelling attempts, yet choosing perpendicular to the gradient works in this model?  We suggest a few reasons this could be, and we leave the more detailed analysis to a future study.  First, there could be some preferential slippage during impact in the perpendicular direction, however if there is any it is likely to be small and we do not foresee it as a probable reason for this particular system.  As a more likely candidate, we could instantaneously have a force opposite the gradient, but on aggregate after several impacts the droplet position ends up being perpendicular to the direction of the wavefield from the previous iteration, which is what seems to happen via inspection in the experiment of S\'{a}enz \ea \cite{Saenz-Mirage2018}.  Finally, the instantaneous wavefield is quite irregular, as S\'{a}enz \ea mention in their article \cite{Saenz-Mirage2018}, and only after averaging the wavefield over many periods do we observe the Mathieu function structure.  Therefore, the droplet could potentially be moving opposite to the gradient of instantaneous irregular wavefields, which sums up to be perpendicular to the Mathieu functions gradient.

Due to the multiple open questions that this communication brings forth, there are many directions for future research.  The most obvious being the detailed hydrodynamic analysis discussed in the previous paragraph.  Another related question is on the precise shape of the instantaneous wavefields at impacts.  This would be a quite difficult experimental problem requiring one or two high speed cameras and many hours of video data.  One technique that could drastically reduce the cost of this analysis is machine learning.  We can imagine capturing the wavefield using a regular camera, which is then segregated into a training set and test set to learn the eigenmodes that contribute to the wavefield.  While there is no guarantee for it to work, if it does it would be significantly less expensive than the full fledged experiment.


\section*{Acknowledgements}
The problem in the present communication was first conceived by A. R., A. U. Oza, and D. L. Blackmore.  In addition to fruitful discussions about the problem, A. R. is indebted to his late advisor, D. L. Blackmore, for his mentorship and support.  A. R. is also grateful to A. U. Oza for the initial progress and discussions about the problem.  The authors also thank T. Rosato for the invitation to contribute to the special issue in honor or D. L. Blackmore.  Finally, A. R. and G. F. Q. appreciate the support of the Amath department at UW, and G. F. Q. appreciates the support of the Physics department at UW.

\bibliographystyle{unsrt}
\bibliography{Bouncing_droplets}

\end{document}